\documentclass[final]{siamltex}

\usepackage{latexsym, amsmath, amssymb}
\usepackage{graphicx,epsfig}
\newtheorem{thm}{Theorem}[section]

\newtheorem{lem}[thm]{Lemma}

\numberwithin{equation}{section}
\newcounter{saveeqn}

\newcommand{\eqnref}[1]{(\ref {#1})}
\newcommand{\beq}{\begin{equation}}
\newcommand{\eeq}{\end{equation}}

\newcommand{\p}{\partial}
\newcommand{\la}{\langle}
\newcommand{\ra}{\rangle}

\newcommand{\qed}{\hfill $\Box$ \medskip}
\newcommand{\RR}{\mathbb{R}}
\newcommand{\Ccal}{\mathcal{C}}
\newcommand{\Dcal}{\mathcal{D}}
\newcommand{\Lcal}{\mathcal{L}}
\newcommand{\Rcal}{\mathcal{R}}

\def\Be{{\bf e}}

\def\Bj{{\bf j}}

\def\Bn{{\bf n}}

\def\Bp{{\bf p}}

\def\Bt{{\bf t}}

\def\Bv{{\bf v}}

\def\Bx{{\bf x}}

\def\BA{{\bf A}}

\def\BD{{\bf D}}

\def\BI{{\bf I}}
\def\BJ{{\bf J}}

\def\BM{{\bf M}}

\def\BO{{\bf O}}
\def\BP{{\bf P}}

\def\BR{{\bf R}}

\def\BT{{\bf T}}

\newcommand{\Ga}{\alpha}

\newcommand{\Gf}{\phi}

\newcommand{\Gc}{\chi}

\newcommand{\Gm}{\mu}

\newcommand{\Gt}{\theta}

\newcommand{\Gs}{\sigma}

\newcommand{\Gy}{\psi}

\newcommand{\GO}{\Omega}

\begin{document}

\title{Bounds on the size of an inclusion using the translation method for two-dimensional complex conductivity\thanks{\footnotesize This work is supported by the Korean Ministry of Education, Sciences and Technology through NRF grants Nos. 2010-0004091 and 2010-0017532, BK21+ at Inha University, and by the National Science Foundation through grant DMS-1211359}}

\author{Hyeonbae Kang\thanks{Department of Mathematics, Inha University, Incheon
402-751, Korea (hbkang, kskim, hdlee@inha.ac.kr, xiaofeili@inha.edu).} \and Kyoungsun Kim\footnotemark[2] \and
Hyundae Lee\footnotemark[2] \and Xiaofei Li\footnotemark[2] \and Graeme W.
Milton\thanks{\footnotesize Department of Mathematics, University
of Utah, Salt Lake City, UT 84112, USA (milton@math.utah.edu).}}

\maketitle

\begin{abstract}
The size estimation problem in electrical impedance tomography is considered when the conductivity is a complex number and the body is two-dimensional. Upper and lower bounds on the volume fraction of the unknown inclusion embedded in the body are derived in terms of two pairs of voltage and current data measured on the boundary of the body. These bounds are derived using the translation method. We also provide numerical examples to show that these bounds are quite tight and stable under measurement noise. 
\end{abstract}

\noindent{\footnotesize {\bf AMS subject classifications}. 65N21, 35J20, 45Q05


\section{Introduction}

The size estimation problem in electrical impedance tomography (EIT) is to estimate the size (area or volume) of unknown inclusions embedded in a conducting body by means of boundary measurements of the voltage and current. The unknown inclusions may represent anomalies in EIT imaging or non-destructive testing or a phase in two phase composite materials. Here we consider the problem where the body is two-dimensional.

To put the problem in a precise way, let $\GO$ be a body in $\RR^2$ occupied by a conducting material and let $D$ be a conducting inclusion inside $\GO$. Let $\Gs_1$ and $\Gs_2$ ($\Gs_1 \neq \Gs_2$) be the conductivities (or dielectric constants) of $D$ and $\GO\setminus D$, respectively, and $\Gs$ be the conductivity profile of $\GO$, {\it i.e.},
 \beq\label{condprof1}
 \Gs= \Gs_1 \Gc(D) + \Gs_2 \Gc(\GO\setminus \overline D)
 \eeq
where $\chi(D)$ is the characteristic function of $D$. If $\GO$ is a two phase composite, we may write $\Gs$ as
\beq\label{condprof2}
\Gs= \Gs_1 \chi_1 + \Gs_2 \chi_2
\eeq
where $\chi_1=1$ in phase 1 and $0$ in phase 2, and $\chi_2=1-\chi_1$. We consider the boundary value problem of the conductivity equation assuming that the Dirichlet data $\phi$ is assigned on $\p\GO$. So the problem to be considered is
\beq\label{PDE}
\begin{cases}
\nabla \cdot \Gs \nabla u = 0 \quad &\mbox{in  } \GO,\\
u=\Gf \ &\mbox{on } \p \GO.
\end{cases}
\eeq
Then the current
\beq
q:= \Gs \nabla u \cdot \Bn
\eeq
is measured on $\p \GO$ where $\Bn$ is the unit outward normal to $\p \GO$. Then the size estimation problem is to estimate the area or volume $|D|$ of the inclusion (or the volume fraction) in terms of a single or finitely many pairs of Cauchy data $(\phi, q)$. It is worth mentioning that we may apply a current on the boundary and measure the corresponding voltage, and methods developed in this this paper can be applied to such situation.

There has been some significant work on the size estimation problem in the context of the conductivity equation. Upper and lower bounds of $|D|$ were obtained by Kang-Seo-Sheen \cite{KSS97}, Alessandrini-Rosset \cite{AR98}, and Alessandrini-Rosset-Seo \cite{ARS00}. These bounds were obtained using estimates of elliptic partial differential equations and expressed by integrals evaluated by a single pair of Cauchy data. A different kind of bound was obtained by Capdeboscq-Vogelius \cite{CV03} using variational methods. Their bounds hold asymptotically when $|D|$ is small. They require special Cauchy data. For special Cauchy data, such as affine boundary conditions on the potential,
the universal bounds of Nemat-Nasser and Hori \cite{nemat} may be inverted to bound $|D|$. Milton \cite{Milton12}, generalizing the results of Nemat-Nasser and Hori, showed that bounds on the properties of composites imply bounds on the response of bodies with special Cauchy data, and these too can be inverted to bound $|D|$ and do not require
the assumption that $|D|$ is small.

Recently, a completely different method to derive bounds on the volume fraction has been introduced by Kang, Kim and Milton which uses translations of the classical variational principles. The translation method was introduced by Murat-Tartar \cite{mt1, mt2, mt3} and Lurie-Cherkaev \cite{lc1, lc2}, and has been used in an essential way to derive bounds on the effective properties of two phase composites in terms of the volume fraction. It turns out that this method of translation can be applied effectively to derive bounds of the volume fraction in terms of boundary measurements: see Kang-Kim-Milton  \cite{KKM12}, and Kang-Milton \cite{KM13}. Numerical implementations of the bounds presented in \cite{KKM12} show that these bounds work quite well to estimate the volume fraction. These bounds are sharp in the sense that for some geometries and for some boundary data the bounds are attained.

In this paper we deal with the case when the conductivity is a complex number. The subject of EIT imaging using complex conductivity has attracted much attention lately since the imaginary part of the complex conductivity changes depending on frequency and images at different frequencies can be used to generate images of high resolution. We refer to \cite{SW13} and references therein for this direction of research. Our purpose is to derive bounds for the volume fraction using boundary measurements when the conductivity is a complex number.

The derivation of bounds in this paper is based on the variational principle of Cherkaev and Gibiansky \cite{CG94} and the translation method. Let $u$ be the solution to \eqnref{PDE} when $\Gs$ is complex. Then the corresponding electric and current fields are given by
\beq
\Be= -\nabla u :=\Be'+i \Be''
\eeq
and
\beq
\Bj= - \Gs \nabla u :=\Bj'+i \Bj''.
\eeq
In above mentioned paper, a minimizing variational principle is obtained for the field $\begin{bmatrix}\Bj'\\ \Be''\end{bmatrix}$.
We may apply the translation method for this field to derive upper and lower bounds for the volume fraction using two Dirichlet boundary data $\phi_1$ and $\phi_2$. But the bounds obtained in this way depends on the choice of boundary data, and it is necessary to consider measurements corresponding to the boundary data $e^{\Gt_1} \phi_1$ and $e^{\Gt_2} \phi_2$ for all $\Gt_1$ and $\Gt_2$. So, we use the parameterized version of the Cherkaev-Gibiansky variational principle which was obtained by Milton-Seppecher-Bouchitte \cite{msb}. Using this variational principle (and translation) a set of bounds parameterized by $\Gt_1$ and $\Gt_2$ is obtained, and by minimizing (or maximizing) them over $\Gt_1$ and $\Gt_2$ we obtain tighter bounds. We emphasize that only the boundary measurements corresponding to one set of boundary data $\phi_1$ and $\phi_2$ are used to compute the bounds. We perform numerical experiments using the bounds obtained in this paper. Results show that the bounds are quite tight and stable in presence of measurement noise.

There is already some work on size estimation for complex conductivity, both for two and three dimensional bodies, but only using a single pair of Cauchy data unlike the two pairs
we use here. Beretta-Francini-Vessella \cite{BFV11} obtained bounds on the size of the inclusion using elliptic estimates, and Thaler-Milton \cite{TM13} developed a comprehensive set of sharp bounds on the volume fraction based on the splitting method. The splitting method, like the translation method, uses the fact that certain integrals (null-Lagrangians) are known in terms of boundary values, but unlike the translation method does not use variational principles, but instead uses the positivity of the norm of certain fields. It was first used in \cite{ml} in the context of elasticity (see also \cite{KMW13}).

This paper is organized as follows. In section \ref{secVP} we review the variational principle of Cherkaev and Gibiansky and its parameterized version, and introduce some null Lagrangians. In section \ref{secTrans} we use null Lagrangians with parameters to translate the variational principle and compute the minimum. In section \ref{sec:tb} parameters are determined and upper and lower bounds for the volume fraction are derived. Section \ref{sec:numeric} presents results of numerical experiments which show that bounds can be quite tight. We finish the paper with a short conclusion. The appendix is to prove a lemma used in the text.

\section{Variational principle and null-Lagrangian}\label{secVP}

We suppose the conductivity $\Gs$ given in \eqnref{condprof1} or \eqnref{condprof2} is a complex constant of the following form:
\beq
\Gs=\Gs'+ i\Gs'', \quad \Gs_1=\Gs_1'+ i\Gs_1'', \quad \Gs_2=\Gs_2'+ i\Gs_2''.
\eeq
Assume that
 \beq \label{condition_sigma}
 \Gs'>0,\quad |\Gs_1| \ne |\Gs_2|,\quad \Gs_1/\Gs_2 \notin \mathbb{R}.
 \eeq
The second condition in the above is required to guarantee that all four points $(\Gs_1',\Gs_1'')$, $(-\Gs_1',-\Gs_1'')$, $(\Gs_2',\Gs_2'')$, $(-\Gs_2',-\Gs_2'')$ are not on the same circle, and the last condition guarantees that all four points are not on the same straight line. (See Section \ref{sec:tb}.)

Let $u$ be the solution $u$ to \eqnref{PDE} and let
\beq
\Be=-\nabla u \quad \mbox{and}\quad \Bj=-\Gs \nabla u.
\eeq
Then we have
\beq\label{govern}
\nabla \cdot \Bj=0 \quad\mbox{and}\quad \nabla \times \Be=0 \quad \mbox{in  } \GO,
\eeq
and
\beq\label{state}
\Bj=\Gs \Be.
\eeq
Let
 \beq
 \Bj=\Bj' + i \Bj'' \quad \mbox{and}\quad \Be=\Be'+i\Be''.
 \eeq
Then, \eqnref{state} is equivalent to the system of equations
\beq
\begin{cases}
\Bj'= \Gs' \Be' - \Gs'' \Be'',\\
\Bj''= \Gs'' \Be' + \Gs' \Be'',
\end{cases}
\eeq
which is in turn equivalent to the following matrix equation:
\beq\label{relation}
\begin{bmatrix}\Be' \\ \Bj''\end{bmatrix} =\frac{1}{\Gs'}\begin{bmatrix} \BI & \Gs''\BI\\
\Gs''\BI & (\Gs'^2+\Gs''^2)\BI\end{bmatrix}\begin{bmatrix}\Bj'\\ \Be''\end{bmatrix} =:\BD_{JE}\begin{bmatrix}\Bj'\\ \Be''\end{bmatrix},
\eeq
where $\BI$ is $2 \times 2$ identity matrix. We then have the variational principle of Cherkaev and Gibiansky \cite{CG94}: for a given Cauchy datum $(\Gf,q)$ on $\p \GO$,
 \beq\label{CherGibi}
 \Big\la \begin{bmatrix}\Bj'\\ \Be''\end{bmatrix}\cdot \BD_{JE}
 \begin{bmatrix}\Bj'\\ \Be''\end{bmatrix}\Big\ra
 =\min_{ {\underline{\Be''}=-\nabla \underline{u''}\atop \underline{u''}=\Gf'' \textrm{ on } \p \GO}
 \ {\nabla \cdot \underline{\Bj'}=0 \atop \underline{\Bj'}\cdot \Bn=-q' \textrm{ on } \p \GO} }
 \Big\la \begin{bmatrix}\underline{\Bj'}\\ \underline{\Be''}\end{bmatrix}\cdot \BD_{JE}\begin{bmatrix}\underline{\Bj'}\\ \underline{\Be''}\end{bmatrix}\Big\ra,
 \eeq
where $\la f \ra$ denotes the average of $f$ over $\GO$, namely,
$$
\la f\ra:=\frac{1}{|\GO|}\int_{\GO}f.
$$
From now on, we put $\BD_{JE}=\BD$ for ease of notation.

We now introduce a parameter $\Gt$ ($0 \le \Gt < 2\pi$) and parameterized variational principle following \cite{msb}. Let $\widetilde\Bj := e^{i\theta} \Bj$ and $\widetilde\Be := e^{i \theta} \Be$. Then we have
\begin{equation*}
\widetilde\Be = - \nabla (u e^{i \theta})\quad \mbox{and} \quad \widetilde\Bj = \Gs \widetilde\Be,
\end{equation*}
and hence
\begin{equation*}
\nabla \cdot \widetilde\Bj = 0 \quad \mbox{and} \quad \nabla \times \widetilde\Be = 0 \quad \mbox{in } \GO.
\end{equation*}
Then we have
\beq
\begin{bmatrix}\widetilde\Be' \\ \widetilde\Bj''\end{bmatrix} = \BD_{JE}\begin{bmatrix}\widetilde\Bj'\\ \widetilde\Be''\end{bmatrix},
\eeq
and the variational principle:
\beq\label{para-CherGibi}
 \Big\la \begin{bmatrix}\widetilde\Bj'\\ \widetilde\Be''\end{bmatrix}\cdot \BD_{JE}
 \begin{bmatrix}\widetilde\Bj'\\ \widetilde\Be''\end{bmatrix}\Big\ra
 =\min \Big\la \begin{bmatrix}\underline{\widetilde\Bj'}\\ \underline{\widetilde\Be''}\end{bmatrix}\cdot \BD_{JE}\begin{bmatrix}\underline{\widetilde\Bj'}\\ \underline{\widetilde\Be''}\end{bmatrix}\Big\ra,
 \eeq
where the minimization is over the trial fields $\underline{\widetilde\Be''}$ and $\underline{\widetilde\Bj'}$ such that
\beq
\begin{cases}
\underline{\widetilde\Be''} =  - \nabla (\Im (\underline{u} e^{i \theta})),
 \quad \nabla \cdot \underline{\widetilde\Bj'} = 0 \quad \mbox{in } \GO,\\
 \underline{u'}=\Gf',\ \underline{u''}=\Gf'', \quad \underline{\widetilde\Bj'} \cdot \Bn =  - \Re(q e^{i \theta}) \quad \mbox{on } \p \GO.
\end{cases}
\eeq
Here and throughout this paper $\Re(z)$ and $\Im(z)$ stand for the real and imaginary parts of $z$, respectively.

Let $\Gf_j$ ($j=1,2$) be given functions (Dirichlet data) defined on $\p \GO$ and $u_j$ be the solution to \eqnref{PDE} when $\Gf=\Gf_j$. Let $q_j=\Gs \nabla u_j \cdot \Bn|_{\p \GO}$. Then $\Be_j=-\nabla u_j$ and $\Bj_j=-\Gs\nabla u_j$ satisfy \eqnref{govern}  and \eqnref{state}, and  $\Bj_j\cdot \Bn=-q_j$ on $\p \GO$. Set
 \beq
\Bv_j = \begin{bmatrix} \widetilde\Bj_j' \\ \widetilde\Be_j'' \end{bmatrix},\quad j=1,2,
\eeq
where $\widetilde\Bj_j := e^{i\theta_j} \Bj$ and $\widetilde\Be_j := e^{i \theta_j} \Be_j$.
The measurement (response) matrix is given by $\BA=(a_{jk})_{j,k =1, 2}$ where
 \beq\label{a}
 a_{jk}=\la \Bv_j\cdot \BD \Bv_k\ra.
 \eeq
We emphasize that $a_{jk}$ is a null-Lagrangian, {\it i.e.}, it can be computed from the boundary measurements. In fact, we have from \eqnref{relation} that
 \begin{align}
a_{jk} & = \frac{1}{|\GO|} \int_{\GO} \begin{bmatrix}\widetilde\Bj_j' \\ \widetilde\Be_j'' \end{bmatrix} \cdot \BD \begin{bmatrix}\widetilde\Bj_k' \\ \widetilde\Be_k'' \end{bmatrix}= \frac{1}{|\GO|}\int_{\GO} \begin{bmatrix}\widetilde\Bj_j' \\ \widetilde\Be_j'' \end{bmatrix} \cdot \begin{bmatrix}\widetilde\Be_k' \\ \widetilde\Bj_k'' \end{bmatrix} \nonumber \\
& = \frac{1}{|\GO|}\int_{\p \GO}\Big[\Re(q_j e^{i \theta_j}) \Re(\Gf_k e^{i \theta_k}) + \Im(q_k e^{i \theta_k}) \Im (\Gf_j e^{i \theta_j})\Big] ds.
\label{measurement}
\end{align}
It is worth mentioning that the measurement matrix $\BA$ depends on the two independent parameters $\Gt_1$ and $\Gt_2$.

Let
 \begin{equation*}
 \BR_{\perp}=\begin{bmatrix} 0& 1\\ -1 & 0\end{bmatrix},
 \end{equation*}
and define for real numbers $t_1$ and $t_2$
 \beq
 \Rcal = \Rcal(t_1, t_2) :=\begin{bmatrix} t_1 \BR_{\perp} & \BO \\ \BO & t_2 \BR_{\perp} \end{bmatrix}.
 \eeq
Since $\Rcal^T=-\Rcal$, we have $\la \Bv_j \cdot \Rcal \Bv_j \ra =0$ for $j=1,2$. Let
 \beq\label{defineb}
 b = b(t_1, t_2):= \la \Bv_1 \cdot \Rcal \Bv_2 \ra = - \la \Bv_2 \cdot \Rcal \Bv_1 \ra.
 \eeq
Then $b$ can be written as
 \beq
 b= \Ga_1 t_1 + \Ga_2 t_2,
 \eeq
where
\beq
\Ga_1 := \la \widetilde\Bj_1' \cdot \BR_{\perp} \widetilde\Bj_2' \ra,\quad \Ga_2 := \la \widetilde\Be_1'' \cdot \BR_{\perp} \widetilde\Be_2'' \ra.
\eeq
We emphasize that $\Ga_1$ and $\Ga_2$ can be computed using the boundary data. In fact, since $\nabla \times \BR_{\perp}\Bj_2 = -\nabla \cdot \Bj_2=0$, there is a potential $\Gy_2$ such that $\BR_{\perp}\Bj_2=\nabla \Gy_2$ in $\GO$. Thus, if $\Bt$ denotes the unit tangent vector on $\p\GO$, then we have
\beq
\Bt\cdot \nabla \Gy_2 = \Bt \cdot \BR_{\perp} \Bj_2=-\Bn \cdot \Bj_2=q_2 \quad \mbox{on } \p\GO.
\eeq
So the boundary value $\Gy_2^0$ of $\Gy_2$ on $\p \GO$ is given by
$$
\Gy_2^0(\Bx)=\int_{\Bx_0}^{\Bx} q_2 ds, \quad \Bx\in\p \GO,
$$
where the integration is along $\p\GO$ in the positive orientation (counterclockwise). Hence
\begin{align}
\Ga_1 & = \frac{1}{|\GO|} \int_{\GO} \widetilde\Bj_1' \cdot (\cos\theta_2 \nabla \Gy_2' - \sin \theta_2 \nabla \Gy_2'') d \Bx \nonumber \\
& = \frac{1}{|\GO|}\int_{\p \GO} \widetilde\Bj_1' \cdot \Bn (\cos \theta_2 \Gy_2' - \sin \theta_2 \Gy_2'') ds (\Bx) \nonumber \\
& =  - \frac{1}{|\GO|}\int_{\p \GO}\Re(q_1 e^{i \theta_1})\Re(\Gy_2^0 e^{i \theta_2})ds (\Bx) \nonumber \\
& = - \frac{1}{|\GO|}\int_{\p \GO} \Re(q_1 e^{i \theta_1})\Re(\int_{\Bx_0}^{\Bx} q_2 e^{i \theta_2} \, ds) ds(\Bx).
\end{align}
Since $\nabla \cdot (\BR_{\perp}\nabla u_2)=0$, we also have
\begin{align}
\Ga_2 & = \frac{1}{|\GO|}\int_{\GO}(\cos \theta_1 \nabla u_1'' + \sin \theta_1 \nabla u_1')\cdot (\cos \theta_2 \BR_{\perp} \nabla u_2'' + \sin \theta_2 \BR_{\perp} \nabla u_2') \nonumber \\
& = \frac{1}{|\GO|}\int_{\p \GO}(u_1'' \cos \theta_1 + u_1' \sin \theta_1)\frac{\p}{\p \Bt}(u_2'' \cos \theta_2 + u_2' \sin \theta_2 )\, ds \nonumber \\
& = \frac{1}{|\GO|}\int_{\p \GO} \Im (\Gf_1 e^{i \theta_1}) \Im \Big(\frac{\p \Gf_2}{\p \Bt} e^{i \theta_2}\Big)\, ds.
\end{align}

\section{Translation of the variational principle}\label{secTrans}

We now apply the translation method to derive bounds for $f_1$, the volume fraction of the phase 1.

We first note that $\begin{bmatrix} 0 & t_3 \BI\\ t_3 \BI & 0 \end{bmatrix}$ applied to fields $\begin{bmatrix}\widetilde\Bj'\\ \widetilde\Be''\end{bmatrix}$ is a null Lagrangian for any real number $t_3$. Define $\Lcal = \Lcal(t_1, t_2, t_3)$ to be the translation of $\begin{bmatrix} \BD & 0 \\ 0 & \BD \end{bmatrix}$ by a null Lagrangian:
 \beq
 \Lcal :=\begin{bmatrix} \widetilde\BD & \Rcal \\ -\Rcal & \widetilde\BD \end{bmatrix}
 \eeq
where
\beq
\widetilde{\BD}:=\BD + \begin{bmatrix} 0 & t_3 \BI\\ t_3 \BI & 0 \end{bmatrix}.
\eeq
We only consider parameters $t_1,~t_2,~t_3$ for which $\Lcal$ is positive semi-definite.

Let
\beq\label{W}
W :=\Big\la \begin{bmatrix}k_1 \Bv_1+k_2 \Bv_2\\ k_3 \Bv_1+k_4 \Bv_2\end{bmatrix} \cdot \Lcal \begin{bmatrix} k_1 \Bv_1+k_2 \Bv_2 \\ k_3\Bv_1+k_4\Bv_2\end{bmatrix}\Big\ra
\eeq
for real numbers $k_1, \ldots, k_4$. One can see that
\beq\label{WDcal}
W = \begin{bmatrix} k_1\\ \vdots \\ k_4\end{bmatrix}\cdot \begin{bmatrix} \widetilde\BA & b \BR_{\perp}\\ -b \BR_{\perp} & \widetilde\BA \end{bmatrix} \begin{bmatrix} k_1\\ \vdots \\ k_4\end{bmatrix},
\eeq
where $\tilde{\BA}=(\tilde{a}_{jk})$ with
\beq\label{definea}
\tilde{a}_{jk}:=a_{jk}+t_3\la \widetilde\Bj_j'\cdot \widetilde\Be_k'' + \widetilde\Be_j'' \cdot \widetilde\Bj_k' \ra,
\eeq
and $b$ is the number defined by \eqnref{defineb}. We emphasize that the new quantity $\tilde{a}_{jk}$ is also determined by the boundary measurements since
\begin{align}
\la \widetilde\Bj_j'\cdot \widetilde\Be_k'' + \widetilde\Be_j'' \cdot \widetilde\Bj_k' \ra & = \frac{1}{|\GO|}\int_{\GO}\Big[- \widetilde \Bj_j'  \cdot \nabla \Im (u_k e^{i \theta_k}) - \nabla \Im(u_j e^{i \theta_j}) \cdot \widetilde \Bj_k'\Big] d \Bx\nonumber\\
& = - \frac{1}{|\GO|}\int_{\p \GO} \Big[(\widetilde\Bj_j' \cdot \Bn) \Im (u_k e^{i \theta_k}) + (\widetilde\Bj_k' \cdot \Bn)  \Im (u_j e^{i \theta_j})\Big] ds\nonumber\\
& = \frac{1}{|\GO|}\int_{\p \GO} \Big[\Re (q_j e^{i \theta_j}) \Im (\Gf_k e^{i \theta_k}) + \Re (q_k e^{i \theta_k}) \Im (\Gf_j e^{i \theta_j})\Big] ds.
\end{align}
Let
 \beq\label{Dcalform}
 \Dcal = \begin{bmatrix} \widetilde\BA & b \BR_{\perp}\\ -b \BR_{\perp} & \widetilde\BA \end{bmatrix}.
 \eeq

Since null Lagrangians are determined by boundary values, one can see from \eqnref{CherGibi} that the following variational principle holds:
\begin{align}
W = \min
 \Big\la\begin{bmatrix}k_1 \underline{\Bv_1}+k_2 \underline{\Bv_2}\\ k_3 \underline{\Bv_1}+k_4 \underline{\Bv_2}\end{bmatrix}\cdot \Lcal \begin{bmatrix}k_1 \underline{\Bv_1}+k_2 \underline{\Bv_2} \\ k_3 \underline{\Bv_1}+k_4 \underline{\Bv_2}\end{bmatrix}\Big\ra,
\end{align}
where the minimum is taken over all $\underline{\Bv_j} =\begin{bmatrix} \underline{\widetilde\Bj_j'} \\ \underline{\widetilde\Be_j''} \end{bmatrix}$, $j=1,2$, satisfying
\begin{equation}\label{constraint}\begin{cases}
\underline{\widetilde\Be''_j} = - \nabla \Im (\underline{u_j} e^{i \theta_j}),
 \quad \nabla \cdot \underline{\widetilde\Bj'_j} = 0 \quad \mbox{in } \GO,\\
 \underline{u'_j}=\Gf_j',\ \underline{u_j''}=\Gf_j'', \quad \underline{\widetilde\Bj_j'} \cdot \Bn = - \Re(q_j e^{i \theta_j}) \quad \mbox{on } \p \GO.\end{cases}
\end{equation}
If $\widetilde\Be_j''$ and $\widetilde\Bj_j'$ satisfy \eqnref{constraint}, one can see that
\begin{align*}
\la \underline{\widetilde\Bj_j'} \ra & = \frac{1}{|\GO|}\int_{\p \GO} \Bx \underline{\widetilde\Bj_j'}\cdot \Bn ds = - \frac{1}{|\GO|}\int_{\p \GO} \Bx \Re(q_j e^{i \theta_j}) = \la \widetilde\Bj_j' \ra,\\
\la \underline{\widetilde\Be_j''}\ra & = - \frac{1}{|\GO|}\int_{\p \GO} \Im (\underline{u_j} e^{i \theta_j}) \Bn = - \frac{1}{|\GO|}\int_{\p \GO} \Im (\Gf_j e^{i \theta_j}) \Bn = \la \widetilde\Be_j'' \ra.
\end{align*}
Hence by relaxing the constraints \eqnref{constraint} for minimization we have
\beq\label{variational}
W \geq \min_{{ \la \underline{\Bv_j}\ra=\la \Bv_j\ra}} \Big\la \begin{bmatrix}
	k_1 \underline{\Bv_1}+k_2 \underline{\Bv_2}   \\
	k_3 \underline{\Bv_1} + k_4 \underline{\Bv_2}
\end{bmatrix} \cdot \Lcal \begin{bmatrix} k_1 \underline{\Bv_1}+k_2 \underline{\Bv_2}\\ k_3 \underline{\Bv_1}+k_4 \underline{\Bv_2}\end{bmatrix}\Big\ra.
\eeq
Here the existence of minimum is guaranteed by the positive semi-definiteness of $\Lcal$.

If the pair $(\hat\Bv_1, \hat\Bv_2)$ is a minimizer of the righthand side of \eqnref{variational}, then we have
\begin{align*}
0 & = \frac{d}{dt}\Big|_{t=0} \Big\la \begin{bmatrix} k_1 (\hat\Bv_1+t\Gy_1)+k_2 (\hat\Bv_2+t\Gy_2) \\ k_3 (\hat\Bv_1+t\Gy_1)+k_4 (\hat\Bv_2+t\Gy_2) \end{bmatrix} \cdot \Lcal \begin{bmatrix} k_1 (\hat\Bv_1+t\Gy_1)+k_2 (\hat\Bv_2+t\Gy_2)\\ k_3 (\hat\Bv_1+t\Gy_1)+k_4 (\hat\Bv_2+t\Gy_2)\end{bmatrix}\Big\ra\\
& = 2\Big\la \begin{bmatrix} k_1 \Gy_1 + k_2 \Gy_2 \\ k_3 \Gy_1 + k_4 \Gy_2 \end{bmatrix} \cdot \Lcal \begin{bmatrix}k_1 \hat\Bv_1+k_2 \hat\Bv_2\\ k_3 \hat\Bv_1+k_4 \hat\Bv_2\end{bmatrix}\Big\ra
\end{align*}
for any pair $(\Gy_1, \Gy_2)$ satisfying $\la \Gy_1\ra = \la \Gy_2 \ra=0$. Thus we have
\beq \label{eq28}
\Lcal \begin{bmatrix}k_1 \hat\Bv_1+k_2 \hat\Bv_2\\ k_3 \hat\Bv_1+k_4 \hat\Bv_2\end{bmatrix}={\bf \Gm}
\eeq
for some constant vector ${\bf \Gm}$.

Let $\Lcal_1$ and $\Lcal_2$ be restrictions of $\Lcal$ to phase 1 and phase 2, respectively, {\it i.e.},
\begin{align*}
\Lcal=\Lcal_1\chi_1 +\Lcal_2 \chi_2.
\end{align*}
Note that $\Lcal_1$ and $\Lcal_2$ are $8 \times 8$ constant matrices. The relation \eqnref{eq28} says that the component of $\begin{bmatrix}k_1 \hat\Bv_1+k_2 \hat\Bv_2\\ k_3 \hat\Bv_1+k_4 \hat\Bv_2\end{bmatrix}\chi_1$ which is orthogonal to $\ker\Lcal_1$ is constant. Likewise, the component of $\begin{bmatrix}k_1 \hat\Bv_1+k_2 \hat\Bv_2\\ k_3 \hat\Bv_1+k_4 \hat\Bv_2\end{bmatrix}\chi_2$ orthogonal to $\ker\Lcal_2$ is constant. Since components in $\ker\Lcal_1$ and $\ker\Lcal_2$ do not contribute to
minimum value in \eqnref{variational}, we obtain
\beq\label{variational45}
W \geq \min_{f_1V_1+f_2V_2=V}
\Big(f_1V_1\cdot\Lcal_1 V_1 + f_2V_2\cdot\Lcal_2 V_2\Big),
\eeq
where $V_1$ and $V_2$ are constant vectors and
\beq
V:=\begin{bmatrix}k_1 \la\Bv_1\ra+k_2 \la\Bv_2\ra\\ k_3 \la\Bv_1\ra+k_4 \la\Bv_2\ra\end{bmatrix}.
\eeq

We use the following lemma whose proof will be given in Appendix.
\begin{lem}\label{lem}
Let $V$ be a finite dimensional vector space, $\Lcal_1,~ \Lcal_2 : V \to V$ self-adjoint linear operators, $f_1$ and $f_2$ positive numbers, and $E_0\in V$. Then
\beq
\min_{f_1E_1 + f_2E_2=E_0}\left( f_1 E_1\cdot \Lcal_1 E_1 + f_2 E_2\cdot \Lcal_2 E_2 \right) =(\pi E_0) \cdot \left[\pi(f_1 \Lcal_1^{-1}+f_2\Lcal_2^{-1})\pi\right]^{-1}\pi E_0,
\eeq
where $\pi$ is the orthogonal projection onto $\mbox{\rm Range}~\Lcal_1 \cap \mbox{\rm Range}~\Lcal_2$ and all the inverses are pseudo-inverses.
\end{lem}

Let $\pi$ be the orthogonal projection onto $\mbox{\rm Range}~\Lcal_1 \cap \mbox{\rm Range}~\Lcal_2$.
Using Lemma \ref{lem}, we know that the minimum on the righthand side of \eqnref{variational45} is
$ V \cdot \Lcal_* V$
where
\beq
\Lcal_*:=\pi \left(\pi(f_1 \Lcal_1^{-1}+f_2 \Lcal_2^{-1})\pi\right)^{-1}\pi .
\eeq
So, we have
$$
W \ge V \cdot \Lcal_* V.
$$
We finally obtain from \eqnref{WDcal} and \eqnref{Dcalform} that
\beq\label{ineq}
\begin{bmatrix} k_1\\ \vdots \\ k_4\end{bmatrix}\cdot \Dcal \begin{bmatrix} k_1\\ \vdots \\ k_4\end{bmatrix}
\geq
 \begin{bmatrix}k_1 \la \Bv_1 \ra + k_2 \la \Bv_2 \ra\\ k_3 \la \Bv_1 \ra + k_4\la \Bv_2 \ra\end{bmatrix} \cdot \Lcal_*  \begin{bmatrix}k_1 \la \Bv_1 \ra + k_2 \la \Bv_2 \ra\\ k_3 \la \Bv_1 \ra + k_4 \la \Bv_2 \ra\end{bmatrix} .
\eeq

We emphasize that $\Lcal_*$ depends on the parameters $t_1, t_2, t_3$. We will choose these parameters in a special way and calculate the corresponding $\Lcal_*$ in the next section. In doing so, the following observation plays a crucial role. Let
\beq
\BJ: = \frac{1}{\sqrt{2}}
\begin{bmatrix}
~1~&0&0&0&0&0&1&0 \\ 0&0&~1~&0&-1&0&0&0 \\ 0&~1~&0&0&0&0&0&1 \\ 0&0&0&~1~&0&-1&0&0 \\
0&0&~1~&0&1&0&0&0 \\ 1&0&0&0&0&0&-1&0 \\ 0&0&0&1&0&1&0&0 \\ 0&1&0&0&0&0&0&-1
\end{bmatrix}.
\eeq
The matrix $\BJ$ has very special properties: it is an orthogonal matrix, namely, $\BJ \BJ^T = \BI$, and the following holds:
\beq\label{L}
\Lcal = \BJ \begin{bmatrix} \BD_{t_3} + \BT & \BO & \BO & \BO\\
\BO & \BD_{t_3} - \BT & \BO & \BO \\
\BO & \BO & \BD_{t_3} + \BT & \BO \\
\BO & \BO & \BO & \BD_{t_3} - \BT \end{bmatrix} \BJ^T \, ,
\eeq
where
\beq
\BD_{t_3} : = \begin{bmatrix}\frac{1}{\Gs'} &  \frac{\Gs''}{\Gs'} + t_3\\  \frac{\Gs''}{\Gs'} + t_3 &  \frac{\Gs'^2 + \Gs''^2}{\Gs'} \end{bmatrix}, \quad
\BT = \begin{bmatrix} t_1 & 0 \\ 0 & t_2 \end{bmatrix}.
\eeq

\section{Translation bounds}\label{sec:tb}
One can see from \eqnref{L} that
$\Lcal \ge 0$ if and only if
\beq\label{posicond}
\BD_{t_3} + \BT \ge 0,\quad \BD_{t_3} - \BT \ge 0.
\eeq
We choose parameters $(t_1, t_2, t_3)$ so that $\Lcal$ is positive semi-definite, more precisely the sum of the ranks of matrices $\Lcal_1$ and $\Lcal_2$ is minimized. Let
\begin{equation*}
\BP_1^\pm := \BD_{t_3}|_{\rm{phase 1}} \pm \BT, \quad
\BP_2^\pm := \BD_{t_3}|_{\rm{phase 2}}  \pm \BT.
\end{equation*}
Then,  $(t_1, t_2, t_3)$ are chosen to be minimizers of
\beq
\min_{t_1, t_2, t_3} \left[ \mbox{rank}~ \BP_1^+ + \mbox{rank} ~\BP_1^- + \mbox{rank}~ \BP_2^+ + \mbox{rank}~ \BP_2^- \right].
\eeq
Such a rank minimizing condition has been used in \cite{LC1, LC2}.

Since $\mbox{rank} ~ \BP_j^\pm \ge 1$, we have
\beq
\min \left[ \mbox{rank}~ \BP_1^+ + \mbox{rank} ~\BP_1^- + \mbox{rank}~ \BP_2^+ + \mbox{rank}~ \BP_2^- \right]=5,
\eeq
and hence there are four possibilities:
\begin{align}
& \det \BP_1^+> 0,\quad \det \BP_1^-=0,\quad \det \BP_2^\pm=0, \label{case1} \\
& \det \BP_2^+> 0,\quad \det \BP_1^\pm=0,\quad \det \BP_2^-=0, \label{case2} \\
& \det \BP_1^-> 0,\quad \det \BP_1^+=0,\quad \det \BP_2^\pm=0, \label{case3} \\
& \det \BP_2^-> 0,\quad \det \BP_1^\pm=0,\quad \det \BP_2^+=0. \label{case4}
\end{align}
The possibilities \eqnref{case1} and \eqnref{case2} yield upper and lower bounds for $f_1$ as we shall see shortly. But, \eqnref{case3} and \eqnref{case4} are equivalent to \eqnref{case1} and \eqnref{case2}, respectively, changing signs of $t_1$ and $t_2$, and hence they yield the same bounds.

Suppose that $(t_1, t_2, t_3)$ satisfies \eqnref{case1}. Following \cite{GM93} (see also \cite[Section 23.7]{mbook}), we interpret this condition in terms of circles.
By explicit calculations, one can see that the last three conditions in \eqnref{case1} are equivalent to the fact that $(-\Gs_1',-\Gs_1''),~(\Gs_2',\Gs_2''),~(-\Gs_2',-\Gs_2'')$ pass through the circle
\beq
t_1(x^2+y^2) + (1+t_1t_2 -t_3^2)x -2t_3 y + t_2 =0.
\eeq
Under the last condition in \eqnref{condition_sigma}, the circle is determined uniquely and  $t_1, ~t_2,~t_3$ are given as follows:
\beq\label{tonettwo}
1/t_1= r \Gs_2'' \pm \sqrt{(r^2+1)|\Gs_2|^2},\quad t_2=-|\Gs_2|^2t_1,\quad t_3= r\Gs_2' t_1
\eeq
where
\beq\label{defr}
r:=\frac{|\Gs_1|^2 - |\Gs_2|^2}{2(\Gs_1'\Gs_2'' - \Gs_2'\Gs_1'')}.
\eeq
Moreover, the second condition in \eqnref{condition_sigma} guarantees the first condition in \eqnref{case1}.

There are additional conditions for $t_1,~t_2, ~t_3$ to fulfill. To ensure \eqnref{posicond}, they should satisfy
\beq
|t_1| \le 1/\Gs_1', \quad |t_1| \le 1/\Gs_2', \quad \det \BP_1^+>0.
\eeq
We show that these conditions can be fulfilled by choosing $t_1$ properly in \eqnref{tonettwo}.

Since $\det \BP_2^+=0$, we have
$$ |\Gs_2|^2( 1/(\Gs_2')^2 - t_1^2 ) = (1/\Gs_2' + t_1 )( |\Gs_2|^2/\Gs_2' + t_2 ) = (\Gs_2'' /\Gs_2' + t_3)^2 \ge 0,$$
and hence $|t_1| \le 1/\Gs_2'$.

On the other hand, since $\det \BP_1^-=0$, we have
\begin{align}
|\Gs_1|^2 \left( \frac{1}{(\Gs_1')^2} - t_1^2 \right) &= \left( \frac{1}{\Gs_1'} - t_1 \right)t_1 (|\Gs_1|^2 - |\Gs_2|^2) +   \left( \frac{1}{\Gs_1'} - t_1 \right)\left( \frac{|\Gs_1|^2}{\Gs_1'} - t_2 \right) \nonumber\\
&=  \left( \frac{1}{\Gs_1'} - t_1 \right)t_1 (|\Gs_1|^2 - |\Gs_2|^2) +   \left( \frac{\Gs_1''}{\Gs_1'} + t_3 \right)^2.\label{eq597}
\end{align}
Let $f$ be a quadratic function whose roots are $1/t_1= r \Gs_2'' \pm \sqrt{(r^2+1)|\Gs_2|^2}$. In fact, it is given by
$$
f(x):=x^2 - 2r \Gs_2'' x - r^2(\Gs_2')^2 - |\Gs_2|^2.
$$
Then one can see that
\beq \label{f_eval}
f(\Gs_1')= - \frac{[((\Gs_1')^2 - (\Gs_2')^2 )|\Gs_2|^2 - (  \Gs_1'\Gs_2'' - \Gs_1'' \Gs_2')^2]^2}{4 ( \Gs_1'\Gs_2'' - \Gs_1'' \Gs_2')^2 (\Gs_2')^2}\le 0.
\eeq
Therefore we have
$$
r \Gs_2'' - \sqrt{(r^2+1)|\Gs_2|^2} \le \Gs_1' \le r \Gs_2'' + \sqrt{(r^2+1)|\Gs_2|^2},
$$
and hence we can choose $t_1$ (among $1/t_1= r \Gs_2'' \pm \sqrt{(r^2+1)|\Gs_2|^2}$) so that
\beq\label{t_crit}
\left( \frac{1}{\Gs_1'} - t_1 \right)t_1 (|\Gs_1|^2 - |\Gs_2|^2)\ge 0 .
\eeq
Then \eqnref{eq597} implies $ |t_1| \le 1/ \Gs_1'$. Here we know $t_1\ne 1/\Gs_1'$ because $t_1=1/\Gs_1'$ would imply $\Gs_2'>\Gs_1'$ so that $0= f(1/t_1)=f(\Gs_1')<0$ by \eqnref{f_eval}. Thus the condition $\det \BP_1^+>0$ (equivalently,  $ t_1(|\Gs_1|^2 - |\Gs_2|^2)>0$) is satisfied automatically with the choice of $t_1$ satisfying \eqnref{t_crit}.

Now we calculate $\Lcal_*$.
First we observe that
\beq\label{nonzero_det}
\det \left ( \BP_1^\pm -\BP_2^\pm\right)=\det \left( \BD_{t_3}|_{\rm{phase 1}} - \BD_{t_3}|_{\rm{phase 2}} \right) = -\frac{(\Gs_1' -\Gs_2')^2 +(\Gs_1''-\Gs_2'')^2}{\Gs_1'\Gs_2'}< 0.
\eeq
Since $\det \BP_1^- = \det \BP_2^- = 0$ while $\BP_1^- - \BP_2^-$ has rank 2, we have
\beq
\mbox{\rm range}~\BP_1^- \cap \mbox{\rm range}~\BP_2^- =0.
\eeq
Since $\det \BP_1^+\ne 0$, we have
\beq
\mbox{\rm range}~\BP_1^+ \cap \mbox{\rm range}~\BP_2^+ =\mbox{\rm range}~\BP_2^+.
\eeq

Recalling from \eqnref{L} that
$$
\Lcal_j = \BJ \begin{bmatrix} \BP_j^+ & \BO & \BO & \BO\\
\BO & \BP_j^- & \BO & \BO \\
\BO & \BO & \BP_j^+ & \BO \\
\BO & \BO & \BO & \BP_j^- \end{bmatrix} \BJ^T \, , \quad j=1,2,
$$
it follows that
\beq
\mbox{\rm range}~\Lcal_1 \cap \mbox{\rm range}~\Lcal_2 =\mbox{\rm range}~\BJ \begin{bmatrix} \BP_2^+ & \BO & \BO & \BO\\
\BO & \BO & \BO & \BO \\
\BO & \BO & \BP_2^+ & \BO \\
\BO & \BO & \BO & \BO \end{bmatrix}\BJ^T \ .
\eeq

Let $\Bp $ be an unit vector generating $\mbox{\rm range} \, \BP_2^+$, and let $\BP$ be the orthogonal projection onto $\mbox{\rm range} \, \BP_2^+$, namely,
$$ \BP= \Bp \Bp^T. $$
Then the orthogonal projection $\pi$ onto $\mbox{\rm range}~\Lcal_1 \cap \mbox{\rm range}~\Lcal_2$  is given by
\beq
\pi=\BJ \begin{bmatrix} \BP & \BO & \BO & \BO\\
\BO & \BO & \BO & \BO \\
\BO & \BO &\BP & \BO \\
\BO & \BO & \BO & \BO \end{bmatrix} \BJ^T .
\eeq
Thus, we have
\begin{align*}
\pi(f_1 \Lcal_1^{-1}+f_2 \Lcal_2^{-1})\pi
&= f_1 \BJ \begin{bmatrix} \BP (\BP_1^+)^{-1} \BP & \BO & \BO & \BO\\
\BO & \BO & \BO & \BO \\
\BO & \BO & \BP (\BP_1^+)^{-1} \BP & \BO \\
\BO & \BO & \BO & \BO \end{bmatrix} \BJ^T \\
& \quad + f_2 \BJ \begin{bmatrix} \BP (\BP_2^+)^{-1} \BP & \BO & \BO & \BO\\
\BO & \BO & \BO & \BO \\
\BO & \BO & \BP (\BP_2^+)^{-1} \BP & \BO \\
\BO & \BO & \BO & \BO \end{bmatrix} \BJ^T \, .
\end{align*}
Here $(\BP_j^+)^{-1}$ is the pseudo-inverse. Since $\BP_2^+$ is symmetric, we have
\beq
\BP_2^+= (\mbox{\rm tr} \BP_2^+) \BP.
\eeq
One can also see that
\beq
\BP (\BP_1^+)^{-1} \BP = \left(\Bp\cdot (\BP_1^+)^{-1}\Bp\right) \BP.
\eeq
Therefore, we have

\begin{align*}
\pi(f_1 \Lcal_1^{-1}+f_2 \Lcal_2^{-1})\pi
&= \left[ f_1 \left(\Bp\cdot (\BP_1^+)^{-1}\Bp\right) + \frac{f_2}{{\rm tr} \BP_2^+} \right] \pi \, ,
\end{align*}
and hence
\beq
\Lcal_*:=\pi \left(\pi(f_1 \Lcal_1^{-1}+f_2 \Lcal_2^{-1})\pi\right)^{-1}\pi =\left[ f_1 \left(\Bp\cdot (\BP_1^+)^{-1}\Bp\right) + \frac{f_2}{{\rm tr} \BP_2^+} \right]^{-1} \pi \, .
\eeq

By positive semi-definiteness of $\BP_1^+$ and $\BP_2^+$, and by \eqnref{case1} and \eqnref{nonzero_det}, we have
\beq\label{posiBp}
\Bp \cdot (\BP_1^+)^{-1}\Bp>0, \quad {\rm tr} \BP_2^+>0.
\eeq
Moreover we have
\begin{align*}
f_1 \left(\Bp\cdot (\BP_1^+)^{-1}\Bp\right) + \frac{f_2}{{\rm tr} \BP_2^+}
&= f_1\left(\Bp \cdot (\BP_1^+)^{-1}\Bp- \frac{1}{{\rm tr} \BP_2^+} \right)+  \frac{1}{{\rm tr} \BP_2^+} \\
 & =-f_1 \frac{\det (\BP_1^+ - \BP_2^+)}{\mbox{\rm tr} \BP_2^+ \det \BP_1^+}+  \frac{1}{{\rm tr} \BP_2^+} .
 \end{align*}
We emphasize that
\beq
-\frac{\det (\BP_1^+ - \BP_2^+)}{\mbox{\rm tr} \BP_2^+ \det \BP_1^+} >0
\eeq
which is a consequence of \eqnref{case1}, \eqnref{nonzero_det} and \eqnref{posiBp}.

Let
\beq
F(f_1):=\left( -f_1 \frac{\det (\BP_1^+ - \BP_2^+)}{\mbox{\rm tr} \BP_2^+ \det \BP_1^+}+  \frac{1}{{\rm tr} \BP_2^+} \right)^{-1}.
\eeq
By \eqnref{ineq} we obtain
\begin{equation}\label{final_positivity} \Dcal=
\begin{bmatrix} \widetilde\BA & b \BR_{\perp}\\ -b \BR_{\perp} & \widetilde\BA \end{bmatrix}  \geq F(f_1)\Ccal^T\begin{bmatrix}\BP & 0 \\ 0 & \BP \end{bmatrix} \Ccal
\end{equation}
where
\beq
\Ccal=\frac{1}{\sqrt{2}} \begin{bmatrix}
\la \widetilde j_{1,1}' \ra & \la \widetilde j_{2,1}' \ra & \la \widetilde j_{1,2}' \ra & \la \widetilde j_{2,2}' \ra \\
\la \widetilde e_{1,1}'' \ra & \la \widetilde e_{2,1}'' \ra & \la \widetilde e_{1,2}'' \ra & \la \widetilde e_{2,2}'' \ra\\
-\la \widetilde j_{1,2}' \ra & -\la \widetilde j_{2,2}' \ra & \la \widetilde j_{1,1}' \ra & \la \widetilde j_{2,1}' \ra\\
-\la \widetilde e_{1,2}'' \ra & -\la \widetilde e_{2,2}'' \ra & \la \widetilde e_{1,1}'' \ra & \la \widetilde e_{2,1}'' \ra
\end{bmatrix} \, .
\eeq
Here $\widetilde j_{k,l}'$ and $\widetilde e_{k,l}''$ ($k=1,2$) are defined by
$$
\widetilde \Bj_{k}' = \begin{bmatrix} \widetilde j_{k,1}' \\ \widetilde j_{k,2}' \end{bmatrix}, \quad \widetilde \Be_{k}'' = \begin{bmatrix} \widetilde e_{k,1}'' \\ \widetilde e_{k,2}'' \end{bmatrix} \, .
$$
We emphasize that $\Ccal$ can be computed using boundary data.

Straightforward calculations show that $\Ccal\begin{bmatrix}\BP & 0 \\ 0 & \BP \end{bmatrix} \Ccal^T$ takes the form
\begin{equation}\label{CPC}
\Ccal^T\begin{bmatrix}\BP & 0 \\ 0 & \BP \end{bmatrix} \Ccal:=\begin{bmatrix} \BM & m \BR_\perp \\ -m \BR_\perp & \BM \end{bmatrix}
\end{equation}
where $\BM$ is a $2 \times 2$ symmetric matrix and $m$ is a real number, which can be computed from boundary values since so does $\Ccal$. Since $\BP$ is singular, we know that $m^2=\det \BM$.
Calculating the eigenvalues of the matrix appearing above, one can see that the inequality \eqnref{final_positivity} is equivalent to the following two inequalities:
\beq\label{fir_ineq}
\mbox{\rm tr} \tilde \BA \ge F(f_1)\mbox{\rm tr} \BM
\eeq
and
\beq\label{sec_ineq}
\det \left(\tilde \BA - F(f_1)\BM\right) \ge (b-F(f_1)m)^2.
\eeq

Inequality \eqnref{fir_ineq} yields a lower bound:
\beq \label{ineq_first}
f_1 \ge -\frac{\mbox{\rm tr} \BP_2^+ \det \BP_1^+}
{\det (\BP_1^+ - \BP_2^+)}\left(\frac{ \mbox{\rm tr}\BM}{\mbox{\rm tr} \tilde \BA } - \frac{1}{\mbox{\rm tr}\BP_2^+ }\right) \, .
\eeq
Note that
\beq
\det \left(\tilde \BA - F(f_1)\BM\right)
	= \det\tilde \BA  - F(f_1)  \mbox{\rm tr} (\tilde\BA \BM^*) + F(f_1)^2 \det \BM,
\eeq
where $\BM^*$ is the adjugate matrix of $\BM$.
So, we obtain from \eqnref{sec_ineq} another lower bound:
\beq \label{ineq_second}
f_1 \ge -\frac{\mbox{\rm tr} \BP_2^+ \det \BP_1^+}
{\det (\BP_1^+ - \BP_2^+)}\left(\frac{\mbox{\rm tr} (\tilde \BA \BM^*)-2bm}{\det \tilde \BA - b^2 } - \frac{1}{\mbox{\rm tr}\BP_2^+ }\right) \, .
\eeq

Observe that $\widetilde\BA, b$ and thus $\BM, m$ depend on $\theta_1$ and $\theta_2$ while $\BP_1$ and $\BP_2$ do not. Denoting the quantities on the righthand sides of inequalities in \eqnref{ineq_first} and \eqnref{ineq_second} by $L_1(\Gt_1, \Gt_2)$ and $L_2(\Gt_1, \Gt_2)$, we have
\beq\label{lowerbound}
f_1 \ge \max_{\theta_1, \theta_2} L_1(\Gt_1, \Gt_2) \vee \max_{\theta_1, \theta_2} L_2(\Gt_1, \Gt_2).
\eeq
Here $a \vee b$ is the maximum of $a$ and $b$. It is worth mentioning that the bound $L_j(\Gt_1, \Gt_2)$ is the same as the bound $L_j(0,0)$ when the boundary data are $e^{\Gt_1} \phi_1$ and $e^{\Gt_2} \phi_2$.

Now suppose that $(t_1, t_2, t_3)$ satisfies \eqnref{case2}. By interchanging the role of phase 1 and phase 2, we obtain
\begin{align}
1-f_1=f_2 &\ge -\frac{\mbox{\rm tr} \BP_1^+ \det \BP_2^+}
{\det (\BP_1^+ - \BP_2^+)}\left( \frac{ \mbox{\rm tr}\BM}{\mbox{\rm tr} \tilde \BA }  - \frac{1}{\mbox{\rm tr}\BP_1^+ }\right),\label{first_upper}\\
1-f_1=f_2 &\ge -\frac{\mbox{\rm tr} \BP_1^+ \det \BP_2^+}
{\det (\BP_1^+ - \BP_2^+)}\left(\frac{\mbox{\rm tr} (\tilde \BA \BM^*)-2bm}{\det \tilde \BA - b^2 } - \frac{1}{\mbox{\rm tr}\BP_1^+ }\right).\label{second_upper}
\end{align}
Here the matrix $\BM$ and the constant $m$ are defined by \eqnref{CPC}, but $\BP$ here is the orthogonal projection onto ${\rm range} \, \BP_1^+$, not ${\rm range} \, \BP_2^+$.

Let
\begin{align*}
U_1(\Gt_1, \Gt_2) &= 1 + \frac{\mbox{\rm tr} \BP_1^+ \det \BP_2^+}
{\det (\BP_1^+ - \BP_2^+)}\left( \frac{ \mbox{\rm tr}\BM}{\mbox{\rm tr} \tilde \BA }  - \frac{1}{\mbox{\rm tr}\BP_1^+ }\right),\\
U_2(\Gt_1, \Gt_2) &= 1 + \frac{\mbox{\rm tr} \BP_1^+ \det \BP_2^+}
{\det (\BP_1^+ - \BP_2^+)}\left(\frac{\mbox{\rm tr} (\tilde \BA \BM^*)-2bm}{\det \tilde \BA - b^2 } - \frac{1}{\mbox{\rm tr}\BP_1^+ }\right).
\end{align*}
Then we have
\beq\label{upperbound}
f_1 \le \min_{\theta_1, \theta_2} U_1(\Gt_1, \Gt_2) \wedge \min_{\theta_1, \theta_2} U_2(\Gt_1, \Gt_2).
\eeq
Here $a \wedge b$ is the minimum of $a$ and $b$.

\section{Numerical experiments}\label{sec:numeric}

This section presents results of some numerical experiments. We compute the bounds for various configurations: (1) the domain is a disk and the inclusion is a concentric disk (Fig. \ref{concentricdisk}, Table \ref{tableconcentricdisk}), (2) domain: a disk, inclusion: an ellipse (Fig. \ref{diskellipse}, Table \ref{tableellipsedisk}), (3) multiple inclusions (Fig. \ref{smile}, Table \ref{tablemultiple}), (4) the domain of general shape (Fig. \ref{generalshape}, Table \ref{tablegeneral}). The results clearly show that bounds obtained in this paper are quite tight, very close to the actual volume fraction. In all computations, we use the Dirichlet boundary data $\phi_1=x$ and $\phi_2=y$, and acquire the corresponding Neumann data by solving \eqnref{PDE} numerically using FEM. We then discretize $[0,2\pi)$ into 200 points, which means $200\times200$ pairs of $(\theta_1,\theta_2)$ are used to optimize the bounds.

We also consider stability of the bounds under measurement noise. In the example of multiple inclusions we add 5, 10, 15, 20\% noise to the Neumann data. We first compute $\nabla u$ by solving \eqnref{PDE} corresponding to the Dirichlet data $\phi_1=x$ and $\phi_2=y$, and then compute
$$
\nabla u^* = [1+ p * \mbox{rand}] \nabla u
$$
for $p=0.05,0.1,0.15, 0.2$ where rand is a generator of Gaussian white noise. So the measured data (with noise) is $q=\nabla u^* \cdot \Bn$. As Table \ref{tablemultiple} shows, the bounds are stable under measurement noise.

Finally we took a configuration (Fig. \ref{disks}) which was considered in \cite{TM13} for the purpose of comparing bounds by the splitting method and those of this paper (translation method). The results presented in Table \ref{tablecomparison} show that the method of this paper yields better bounds than the slitting method. It is worth emphasizing that the splitting method in \cite{TM13} uses a single measurement while the translation method uses two measurements.

\begin{figure}
\centering
 \setlength{\unitlength}{1bp}%
  \begin{picture}(93.03, 93.03)(0,0)
  \put(0,0){\includegraphics{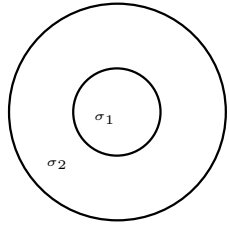}}
  \put(37.78,43.35){\fontsize{4.12}{4.95}\selectfont $\sigma_1$}
  \put(19.85,26.24){\fontsize{4.12}{4.95}\selectfont $\sigma_2$}
  \end{picture}%
\caption{Concentric disks. Inclusion $D$ is the disk of radius $r_1=0.4$ and of conductivity $\Gs_1$. $\Omega$ is the unit disk and the conductivity of $\Omega \setminus \overline D$ is $\sigma_2=1$.}
    \label{concentricdisk}
\end{figure}

\begin{table}[t]
\caption{Two concentric disks. $\sigma_1$: the conductivity of the inclusion $D$, $\sigma_2$: conductivity of $\Omega \setminus \overline D$, $f_1$: the exact area fraction of the inclusion, $max(L)$: the lower bound, $min(U)$: the upper bound.}\label{tableconcentricdisk}
\begin{center}
\begin{tabular}{lllllll}
                 \hline
                 $\sigma_1$ & $\sigma_2$& $f_1$ & $max(L)$ & $max(L)/f_1$ & $min(U)$ & $min(U)/f_1$ \\
                 \hline
1+i	&1&0.16&0.159919&0.999492&0.160044	 &1.000274	\\
2+0.5i	&1&0.16	&0.159944&0.999652	&0.160015	 &1.000097\\
2+5i&1	&0.16  &0.159937&0.999608   &0.160008	 &1.000049	\\
4+100i&1	&0.16	& 0.159839&0.998992&0.160026&1.000165	\\
                 \hline
\end{tabular}
\end{center}
\end{table}

\begin{figure}
\centering
\setlength{\unitlength}{1bp}%
  \begin{picture}(93.55, 93.55)(0,0)
  \put(0,0){\includegraphics{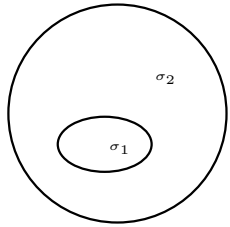}}
  \put(43.73,32.88){\fontsize{4.12}{4.95}\selectfont $\sigma_1$}
  \put(60.97,59.45){\fontsize{4.12}{4.95}\selectfont $\sigma_2$}
  \end{picture}%
\caption{$\Omega$ is the unit disk, D is an ellipse with center point lying at $(-0.1,-0.3)$ with the major axis $a=0.4$ and the minor axis $b=0.3$. The conductivity of inclusion is $\sigma_{1}=2+i$ and the background conductivity is $\sigma_2=1$.}
    \label{diskellipse}
\end{figure}

\begin{table}[t]
\caption{Elliptic inclusion}\label{tableellipsedisk}
\begin{center}
\begin{tabular}{lllllll}
                 \hline
                 $\sigma_1$ & $\sigma_2$& $f_1$ & $max(L)$ & $max(L)/f_1$ & $min(U)$ & $min(U)/f_1$ \\
\hline
2+i&1	&0.12  &0.119559&0.996324& 0.120800	 &1.00667	\\
                 \hline
\end{tabular}
\end{center}
\end{table}

\begin{figure}
\centering
\setlength{\unitlength}{1bp}%
  \begin{picture}(90.99, 90.99)(0,0)
  \put(0,0){\includegraphics{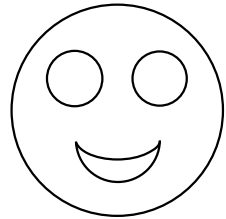}}
  \end{picture}%
  \caption{$\Omega$ is the unit disk, D is composed of three parts: two circles with radius 0.25 centered at (-0.4, 0.3) and (0.4 0.3), and one crescent inclusion with area 0.0225. The conductivity of the three inclusions is $\sigma_1=2+i$ and the background conductivity is $\sigma_2=1$. The exact area is $f_1=0.1475$.}
    \label{smile}
\end{figure}

\begin{table}[t]
\caption{Multiple inclusions. $\sigma_1=2+i$, $\sigma_2=1$, $f_1=0.1475$. }\label{tablemultiple}
\begin{center}
\begin{tabular}{lllll}
                 \hline
                noise level & $max(L)$ & $max(L)/f_1$ & $min(U)$ & $min(U)/f_1$ \\
                 \hline
0  &0.146614& 0.993991 &0.148187 &1.004655\\
5\%& 0.143527&  0.973066 &0.151170 &1.024883\\
10\% &0.134217& 0.909943 &0.159726&1.082891\\
15\% &0.119537& 0.810420 &0.174828&1.185274\\
20\%  &  0.098495&  0.667764 & 0.194967 & 1.321812\\
                 \hline
\end{tabular}
\end{center}
\end{table}

\begin{figure}
\centering
\setlength{\unitlength}{1bp}%
  \begin{picture}(107.72, 78.85)(0,0)
  \put(0,0){\includegraphics{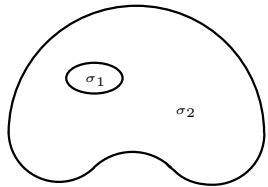}}
  \put(34.56,44.41){\fontsize{4.28}{5.13}\selectfont $\sigma_1$}
  \put(68.55,32.62){\fontsize{4.28}{5.13}\selectfont $\sigma_2$}
  \end{picture}%
  \caption{$\GO$ is of general shape. $D$ is with conductivity $\sigma_{1}=2+i$ and the background conductivity is $\sigma_2=1$. The exact area fraction is $f_1=0.029281$ .}
 \label{generalshape}
 \end{figure}

\begin{table}[t]
\caption{General shape domain}\label{tablegeneral}
\begin{center}
\begin{tabular}{lllllll}
                 \hline
                 $\sigma_1$&$\sigma_2$ & $f_1$ & $max(L)$ & $max(L)/f_1$ & $min(U)$ & $min(U)/f_1$ \\
\hline
1+2i&1&0.029281 &0.029172&0.996278  &0.029631&1.011941	\\
                 \hline
\end{tabular}
\end{center}
\end{table}

\begin{figure}
\centering
\setlength{\unitlength}{1bp}%
  \begin{picture}(99.02, 99.02)(0,0)
  \put(0,0){\includegraphics{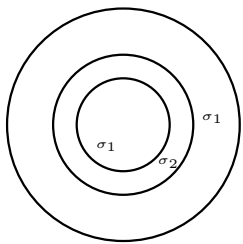}}
  \put(39.27,41.02){\fontsize{5.30}{6.36}\selectfont $\sigma_1$}
  \put(62.58,34.77){\fontsize{5.30}{6.36}\selectfont $\sigma_2$}
  \put(79.07,51.51){\fontsize{5.30}{6.36}\selectfont $\sigma_1$}
  \end{picture}%
  \caption{The configuration from \cite{TM13}. Phase 1 consist of the core and the outer annulus, and its area fraction is $f_1=0.8$. The conductivity of phase 1 is $\sigma_1=3+8i$, and that of phase 2 is $\Gs_2= 8+6i$. The radii of three circles are $R_1=2$, $R_2=3$, $R_3=5$. }
\label{disks}
  \end{figure}

\begin{table}[t]
\caption{Comparison with the splitting method results. $\Gs_1=3+8i$, $\Gs_2=8+6i$, $f_1=0.8$}\label{tablecomparison}
\begin{center}
\begin{tabular}{llllllll}
                 \hline
                 Method  & $max(L)$ & $max(L)/f_1$ & $min(U)$ & $min(U)/f_1$ \\
\hline
splitting method  &0.794&0.9925   &0.808	&1.01	\\
translation method  &0.799485& 0.999356 &0.800064&1.000080	\\
                 \hline
\end{tabular}
\end{center}
\end{table}

\section*{Conclusion}

We have derived upper and lower bounds of the volume fraction of an unknown inclusion (or two phase composites) using boundary measurements when the conductivity is complex. We use the minimizing variational principles with parameters for the fields $\Be= \Be'+i \Be''$ and $\Bj= \Bj'+i \Bj''$. The bounds are given in a nonlinear way in terms of the determinant and the trace of the measurement matrix, and some other null Lagrangians which can be computed using boundary measurements.
We perform numerical experiments to validate the effectiveness of the bounds obtained in this paper and to compare them with those in \cite{TM13}. Results show that the bounds obtained in this paper are quite tight and stable under measurement (white) noise. They also show that these bounds are better than those obtained in \cite{TM13} using less boundary measurement data.

\appendix

\section{Proof of Lemma \ref{lem}}

This appendix is to prove Lemma \ref{lem}.

We consider the following minimization problem:
\beq
\min_{E_1, E_2\in V}\left( f_1 E_1\cdot \Lcal_1 E_1 + f_2 E_2\cdot \Lcal_2 E_2 - 2 A\cdot (f_1E_1 + f_2E_2)\right),
\eeq
where the Lagrange multiplier $A$ is a vector in $\mbox{\rm Range}~\Lcal_1 \cap \mbox{\rm Range}~\Lcal_2$ (otherwise there is no minimum).  If $E_1$ and $E_2$ are minimizers, they should satisfy
\beq 2f_1 \delta E_1 \cdot (\Lcal_1E_1 -A)=0,\quad 2f_2 \delta E_2 \cdot (\Lcal_2E_2 -A)=0
\eeq
for all the increments $\delta E_1$ and $\delta E_2$. Then
\beq E_1= \Lcal_1^{-1}A+ E_1^0,\quad E_2= \Lcal_2^{-1}A+ E_2^0 \eeq
for some $E_1^0\in \ker \Lcal_1$ and  $E_2^0\in \ker \Lcal_2$. Since ${\rm Range}~\Lcal_j$ is orthogonal to $\ker \Lcal_j$, if we impose the constraint $f_1E_1 + f_2E_2=E_0$, then we have
$$
\pi E_0 = \pi (f_1 \Lcal_1^{-1}+f_2\Lcal_2^{-1})\pi A,
$$
and hence
\beq
A=\left[\pi(f_1 \Lcal_1^{-1}+f_2\Lcal_2^{-1})\pi\right]^{-1}\pi E_0.
\eeq
And we have
\begin{align*}
f_1 E_1 \cdot \Lcal_1 E_1 + f_2 E_2 \cdot \Lcal E_2 &= f_1A\cdot \Lcal_1^{-1} A+ f_2 A\cdot \Lcal_2^{-1} A \\
&=(\pi E_0) \cdot \left[\pi(f_1 \Lcal_1^{-1}+f_2\Lcal_2^{-1})\pi\right]^{-1}\pi E_0.
\end{align*}
This completes the proof. \qed


\end{document}